\documentclass[11pt]{article}
\usepackage{amsmath}
\usepackage{amssymb}
\usepackage{amsthm}
\usepackage[usenames]{color}
\usepackage{amscd}
\usepackage[colorlinks=true,
linkcolor=webgreen,
filecolor=webbrown,
citecolor=webgreen]{hyperref}

\definecolor{webgreen}{rgb}{0,.5,0}
\definecolor{webbrown}{rgb}{.6,0,0}

\def\1{{\bf 1}}

\def\N{{\Bbb N}}
\def\Z{{\Bbb Z}}
\def\C{{\Bbb C}}

\def\id{\operatorname{id}}

\newtheorem{theorem}{Theorem}[section]
\newtheorem{corollary}[theorem]{Corollary}

\numberwithin{equation}{section}

\begin{document}

\title{\bf Two generalizations of the Busche-Ramanujan identities}
\author{L\'aszl\'o T\'oth}
\date{}
\maketitle

\centerline{Int. J. Number Theory {\bf 9} (2013), 1301-1311}

\begin{abstract} We derive two new generalizations of the Busche-Ramanujan identities
involving the multiple Dirichlet convolution of arithmetic functions
of several variables. The proofs use formal multiple Dirichlet
series and properties of symmetric polynomials of several variables.
\end{abstract}



\section{Introduction}

Throughout the paper we use the notation: $\N=\{1,2,\ldots\}$,
$\N_0=\{0,1,2,\ldots\}$, $\1(n)=1$, $\id(n)=n$, $\id_k(n)=n^k$
($k\in \C$, $n\in \N$), $\mu$ is the M\"obius function, $*$ denotes
the Dirichlet convolution of arithmetic functions, $f^{-1_*}$ is the inverse under $*$ of the function $f$,
$\zeta$ is the Riemann zeta function.

Let $g$ and $h$ be two completely multiplicative arithmetic
functions and let $f=g*h$. The Busche-Ramanujan identities state that for every
$m,n\in \N$,
\begin{equation} \label{Busche_Raman_i}
f(mn) = \sum_{a\mid \gcd(m,n)} f\left(\frac{m}{a}\right)
f\left(\frac{n}{a}\right)\mu(a)g(a)h(a)
\end{equation}
and
\begin{equation} \label{Busche_Raman_ii}
f(m)f(n) = \sum_{a\mid \gcd(m,n)} f\left(\frac{mn}{a^2}\right)
g(a)h(a).
\end{equation}

For example, these identities hold true for the following special
functions $f$: (i) the function $\sigma_k=\1*\id_k$, in particular,
the divisor function $d=\1*\1$ and the sum-of-divisors function
$\sigma= \1* \id$; (ii) the alternating sum-of-divisors function
$\beta =\lambda * \id$, where $\lambda$ is the Liouville function,
cf. \cite{Tot2012}; (iii) the function $R_1=\1*\chi$, where $\chi$
is the nonprincipal character (mod $4$), $R(n)=4R_1(n)$ representing
the number of ordered pairs $(x,y)\in \Z^2$ such that $n=x^2+y^2$;
(iv) Ramanujan's function $\tau$, defined by the expansion
\begin{equation*}
\sum_{n=1}^{\infty} \tau(n)x^n =x \prod_{n=1}^{\infty} (1-x^n)^{24},
\quad |x|<1,
\end{equation*}
where $\tau = g_1*g_2$ with certain completely multiplicative
functions $g_1, g_2$ verifying $g_1(n)g_2(n)=n^{11}$ for every $n\in
\N$, cf., e.g., \cite{Apo1990,LucShp2003}.

For history and discussion, as well as for generalizations and
analogues of the Busche-Ramanujan identities we refer to
\cite{Bal1997, CarGio1975, Hau1988, Hau1999, Hau1999a, MacWon2012, McC1960, McC1962, McC1986,
Mer1989, RajSit2007, Siv1989, Vai1931}. The functions which are the
convolution of two completely multiplicative functions are called in
the literature specially multiplicative or quadratic functions. See the papers
\cite{Hau2003, RedSiv1981, Siv1976} for various other properties of such functions.
We note that the functions $f$ which can be written as $f=g*h^{-1_*}$, where $g$ and $h$ are
completely multiplicative are called totients. If $f$ is a totient, in particular if $f$ is
the Euler totient function
$\varphi=\id* \mu$, then $f$ verifies a restricted Busche-Ramanujan identity, cf.
\cite{Hau1996, McC1986, Vai1931}.

In the present paper we derive two new generalizations of the
Busche-Ramanujan identities. Namely, we consider the values of a
specially multiplicative function for products of several arbitrary
integers (Theorem \ref{Th_1}), and then deduce formulae for the
convolution of several arbitrary completely multiplicative functions
(Theorem \ref{Th_2}). Our general formulae involve multiplicative
functions of several variables and their Dirichlet convolution. Such
a treatment explains also the equivalence of the identities of type
\eqref{Busche_Raman_i} and \eqref{Busche_Raman_ii}. The proofs use
simple arguments concerning formal multiple Dirichlet series of
arithmetic functions of several variables and properties of
symmetric polynomials of several variables.

%%%%%%%%%%%%%%%%%%%%%%%%%%%%%%%%%%%%%%%%%%%%%%%%%%%%%%%%%%%%%%%%%%%%%%%%%
%%%%%%%%%%%%%%%%%%%%%%%%%%%%%%%%%%%%%%%%%%%%%%%%%%%%%%%%%%%%%%%%%%%%%%%%%

\section{Preliminaries} \label{Sect_Prelim}

In the paper we use the following notions and properties.

\subsection{Arithmetic functions of several variables}

Let ${\cal
F}_r$ be the set of arithmetic functions of $r$ ($r\in \N$)
variables, i.e., of functions $f: \N^r \to \C$. If $f,g\in {\cal
F}_r$, then their Dirichlet convolution is defined by
\begin{equation*}
(f*g)(n_1,\ldots,n_r)= \sum_{a_1\mid n_1, \, \ldots, \, a_r\mid n_r}
f(a_1,\ldots,a_r) g(n_1/a_1, \ldots, n_r/a_r).
\end{equation*}

Here $({\cal F}_r,+,*)$ is a commutative ring with identity
$\delta$, where $\delta(1,\ldots,1)=1$ and
$\delta(n_1,\ldots,n_r)=0$ for $n_1\cdots n_r>1$. A function $f\in
{\cal F}_r$ has an inverse under the convolution, denoted by
$f^{-1_*}$, if and only if $f(1,\ldots,1)\ne 0$.

A function $f\in {\cal F}_r$ is said to be multiplicative if it is
not identically zero and
\begin{equation*}
f(m_1n_1,\ldots,m_rn_r)= f(m_1,\ldots,m_r)f(n_1,\ldots,n_r)
\end{equation*}
holds for any $m_1,\ldots,m_r,n_1,\ldots,n_r\in \N$ such that
$\gcd(m_1\cdots m_r,n_1\cdots n_r)=1$. If $f$ is multiplicative,
then it is determined by the values $f(p^{\nu_1},\ldots,p^{\nu_r})$,
where $p$ is prime and $\nu_1,\ldots,\nu_r\in \N_0$. More exactly,
$f(1,\ldots,1)=1$ and for any $n_1,\ldots,n_r\in \N$,
\begin{equation*}
f(n_1,\ldots,n_r)= \prod_p f(p^{\nu_p(n_1)},\ldots,p^{\nu_p(n_r)}),
\end{equation*}
where $n_i=\prod_p p^{\nu_p(n_i)}$ is the prime power factorization
of $n_i$ ($1\le i\le r$), the products being over the primes $p$ and
all but a finite number of the exponents $\nu_p(n_i)$ are zero.

The convolution of multiplicative functions is multiplicative. The
inverse of a multiplicative function is multiplicative. If $f\in
{\cal F}_1$ is multiplicative, then the function $(n_1,\ldots,n_r)
\mapsto f(n_1\cdots n_r)$ is multiplicative.

The (formal) multiple Dirichlet series of a function $f\in {\cal
F}_r$ is given by
\begin{equation*}
D(f;s_1,\ldots,s_r)= \sum_{n_1,\ldots,n_r=1}^{\infty}
\frac{f(n_1,\ldots,n_r)}{n_1^{s_1}\cdots n_r^{s_r}}.
\end{equation*}

We have
\begin{equation*}
D(f*g;s_1,\ldots,s_r) = D(f;s_1,\ldots,s_r) D(g;s_1,\ldots,s_r),
\end{equation*}
and
\begin{equation*}
D(f^{-1_*};s_1,\ldots,s_r) = D(f;s_1,\ldots,s_r)^{-1},
\end{equation*}
formally or in the case of absolute convergence.

If $f\in {\cal F}_r$ is multiplicative, then its Dirichlet series
can be expanded into a (formal) Euler product, that is,
\begin{equation*}
D(f;s_1,\ldots,s_r)=  \prod_p \sum_{\nu_1,\ldots,\nu_r=0}^{\infty}
\frac{f(p^{\nu_1},\ldots, p^{\nu_r})}{p^{\nu_1s_1+\ldots +\nu_r
s_r}},
\end{equation*}
the product being over the primes $p$.

If $r=1$, i.e., in the case of functions of a single variable we
reobtain the familiar notions and properties, cf.
\cite{Apo1976,McC1986,Siv1989}. See \cite[Ch.\ VII]{Siv1989},
\cite{Tot2011}, \cite{Vai1931} for the case of several variables.

We also need that if $g\in {\cal F}_1$ is completely multiplicative, then
\begin{equation} \label{Dir_ser_g}
D(g;s) = \prod_p \left(1- \frac{g(p)}{p^s} \right)^{-1}.
\end{equation}

\subsection{Symmetric polynomials}

Let $e_d(x_1,\ldots,x_k)=\sum_{1\le i_1<\ldots <i_d\le k}
x_{i_1}\cdots x_{i_d}$ be the elementary symmetric polynomial in
$x_1$, $\ldots,x_k$ of degree $d$ ($d\ge 0$). By convention,
$e_0(x_1,\ldots,x_k)=1$ and $e_d(x_1,\ldots,x_k)=0$ ($d\ge k+1$).
Furthermore, let $h_d(x_1,\ldots,x_k)=\sum_{1\le i_1\le \ldots \le
i_d\le k} x_{i_1}\cdots x_{i_d}$ stand for the complete homogeneous
symmetric polynomial in $x_1,\ldots,x_k$ of degree $d$ ($d\ge 0$).
By convention, $h_0(x_1,\ldots,x_k)=1$.

We will use the polynomial identity (see, e.g., \cite[Ch.\ 8]{Aig2007})
\begin{equation} \label{rel_e_h}
\sum_{d=0}^n (-1)^d e_d(x_1,\ldots,x_k) h_{n-d}(x_1,\ldots,x_k)=
\begin{cases} 1, & n=0, \\ 0, &  n\ge 1. \end{cases}
\end{equation}

We also need the representation
\begin{equation} \label{ident_h}
h_d(x_1,\ldots,x_k)=\sum_{i=1}^k x_i^{d+k-1} \prod_{\substack{j=1\\
j\ne i}}^k (x_i-x_j)^{-1} \quad (d\ge 0),
\end{equation}
cf. \cite{Cor2011}. In the case $k=2$  the identity \eqref{ident_h} reduces to
\begin{equation} \label{ident_h_2}
h_d(x,y)=\frac{x^{d+1}-y^{d+1}}{x-y} \quad (d\ge 0).
\end{equation}

The connection between arithmetic functions and symmetric
polynomials we utilize in the proofs is the following (see, e.g.,
\cite{MacGea2005} for other aspects): If $f_1,\ldots,f_k$ are
completely multiplicative functions of a single variable and
$F=f_1*\cdots
* f_k$, then for every prime power $p^{\nu}$ ($\nu \in \N$),
\begin{equation} \label{F_homogen}
F(p^{\nu})= h_{\nu}(x_1,\ldots,x_k),
\end{equation}
where $x_i=f_i(p)$ ($1\le i\le k$).

%%%%%%%%%%%%%%%%%%%%%%%%%%%%%%%%%%%%%%%%%%%%%%%%%%%%%%%%%%%%%%%%%%%%%%%%%
%%%%%%%%%%%%%%%%%%%%%%%%%%%%%%%%%%%%%%%%%%%%%%%%%%%%%%%%%%%%%%%%%%%%%%%%%

\section{Results}

We prove the following results.

\begin{theorem} \label{Th_1} Let $g$ and $h$ be two completely multiplicative
functions and let $f=g*h$. Let $\psi_f$ be the multiplicative
function of $r$ ($r\in \N$) variables defined as follows. For every
prime $p$ and every $\nu_1,\ldots,\nu_r\in \N_0$ set
\begin{equation*}
\psi_f(p^{\nu_1},\ldots,p^{\nu_r})=
\begin{cases} 1, & \nu_1=\ldots=\nu_r=0,\\
(-1)^{j-1} g(p)h(p)f(p^{j-2}), &
\nu_1,\ldots,\nu_r\in \{0,1\}, \\ & j:=\nu_1+\ldots+\nu_r \ge 2,\\
0, & \text{otherwise}.
\end{cases}
\end{equation*}

Furthermore, let $\psi_f^{-1_*}$  be the inverse under the $r$
variables convolution of the function $\psi_f$. Then for every
$n_1,\ldots,n_r\in \N$ ($r\in \N$),
\begin{equation} \label{Busche_Raman_gener}
f(n_1\cdots n_r) = \sum_{a_1\mid n_1, \, \ldots, \, a_r\mid n_r}
f\left(\frac{n_1}{a_1}\right)\cdots  f\left(\frac{n_r}{a_r}\right)
\psi_f(a_1,\ldots,a_r),
\end{equation}
and
\begin{equation} \label{Busche_Raman_gener_inv}
f(n_1)\cdots f(n_r) = \sum_{a_1\mid n_1,\, \ldots,\, a_r\mid n_r}
f\left(\frac{n_1\cdots n_r}{a_1\cdots a_r} \right)
\psi_f^{-1_*}(a_1,\ldots,a_r).
\end{equation}
\end{theorem}

\begin{theorem} \label{Th_2} Let $f_1,\ldots,f_k$ be completely multiplicative
functions ($k\in \N$) and let $F=f_1*\cdots * f_k$. Let $\vartheta_F$
be the multiplicative function of two variables defined as follows.
For every prime $p$ and every $\nu_1,\nu_2\in \N_0$ set
\begin{equation*}
\vartheta_F(p^{\nu_1},p^{\nu_2})=
\begin{cases} 1, & \nu_1=\nu_2=0,\\
(-1)^{\nu_1+\nu_2-1} e_{\nu_1+\nu_2}(f_1(p),\ldots,f_k(p)), &
\nu_1,\nu_2\ge 1, \nu_1+\nu_2 \le k,\\
0, & \text{otherwise},
\end{cases}
\end{equation*}
where $e_d(x_1,\ldots,x_k)$ represents the elementary symmetric
polynomial in $x_1,\ldots,x_k$ of degree $d$. Furthermore, let
$\vartheta_F^{-1_*}$ denote the inverse under the two variables
convolution of the function $\vartheta_F$. Then for every
$n_1,n_2\in \N$,
\begin{equation} \label{Busche_Raman_gener_other}
F(n_1n_2) = \sum_{a_1\mid n_1,  \, a_2\mid n_2}
F\left(\frac{n_1}{a_1}\right) F \left(\frac{n_2}{a_2}\right)
\vartheta_F(a_1,a_2),
\end{equation}
and
\begin{equation} \label{Busche_Raman_gener_other_inv}
F(n_1) F(n_2) = \sum_{a_1\mid n_1, \, a_2\mid n_2}
F\left(\frac{n_1n_2}{a_1a_2}\right) \vartheta_F^{-1_*}(a_1,a_2).
\end{equation}
\end{theorem}

Theorems \ref{Th_1} and \ref{Th_2} reduce in the cases $r=2$,
respectively $k=2$ to the Busche-Ramanujan identities
\eqref{Busche_Raman_i} and \eqref{Busche_Raman_ii}. Note that if
$r=1$, then Theorem \ref{Th_1} is trivial. It gives $f(n_1)=f(n_1)$
with $\psi_f=\delta$ for every $f$. If $k=1$, then the identities of
Theorem \ref{Th_2} reduce to  $f_1(n_1n_2)=f_1(n_1)f_1(n_2)$, where
$\vartheta_F=\delta$ for every $F=f_1$.

For Ramanujan's tau function and for the Piltz divisor function
$d_k$ we obtain from our results the next identities. We recall that
$d_k(n)$ is the number of ordered $k$-tuples $(x_1,\ldots,x_k)$ of
positive integers such that $x_1\cdots x_k=n$. Therefore,
$d_k=\1*\cdots *\1$ ($k$-times).

\begin{corollary} \label{Cor_1} For every
$n_1,\ldots,n_r\in \N$ ($r\in \N$),
\begin{equation*} \label{Busche_Raman_gener_tau}
\tau(n_1\cdots n_r) = \sum_{a_1\mid n_1, \, \ldots, \, a_r\mid n_r}
\tau\left(\frac{n_1}{a_1}\right)\cdots
\tau\left(\frac{n_r}{a_r}\right) \psi_{\tau}(a_1,\ldots,a_r),
\end{equation*}
where the multiplicative function $\psi_{\tau}$ is defined for every
prime $p$ and every $\nu_1,\ldots,\nu_r\in \N_0$ by
\begin{equation*}
\psi_{\tau}(p^{\nu_1},\ldots,p^{\nu_r})=
\begin{cases} 1, & \nu_1=\ldots=\nu_r=0,\\
(-1)^{j-1} p^{11} \tau(p^{j-2}), &
\nu_1,\ldots,\nu_r\in \{0,1\}, j:=\nu_1+\ldots+\nu_r \ge 2,\\
0, & \text{otherwise}.
\end{cases}
\end{equation*}
\end{corollary}

\begin{corollary} \label{Cor_2} Let $k\in \N$. For every $n_1,n_2\in \N$,
\begin{equation*} \label{Busche_Raman_gener_d_k}
d_k(n_1n_2) = \sum_{a_1\mid n_1,  \, a_2\mid n_2}
d_k\left(\frac{n_1}{a_1}\right) d_k\left(\frac{n_2}{a_2}\right)
\vartheta_k(a_1,a_2),
\end{equation*}
where the multiplicative function $\vartheta_k$ is defined for every
prime $p$ and every $\nu_1,\nu_2\in \N_0$ by
\begin{equation*}
\vartheta_k(p^{\nu_1},p^{\nu_2})=
\begin{cases} 1, & \nu_1=\nu_2=0,\\
(-1)^{\nu_1+\nu_2-1} \binom{k}{\nu_1+\nu_2}, &
\nu_1,\nu_2\ge 1, \nu_1+\nu_2 \le k,\\
0, & \text{otherwise}.
\end{cases}
\end{equation*}
\end{corollary}

%%%%%%%%%%%%%%%%%%%%%%%%%%%%%%%%%%%%%%%%%%%%%%%%%%%%%%%%%%%%%%%%%%%%%%%%%
%%%%%%%%%%%%%%%%%%%%%%%%%%%%%%%%%%%%%%%%%%%%%%%%%%%%%%%%%%%%%%%%%%%%%%%%%

\section{Proof of Theorem \ref{Th_1}} \label{Proof_Theorem_1}

Formula \eqref{Busche_Raman_gener} is a direct consequence of the
following identity concerning multiple Dirichlet series: If $g$ and
$h$ are completely multiplicative functions and $f=g*h$, then
\begin{equation*}
\sum_{n_1,\ldots,n_r=1}^{\infty} \frac{f(n_1\cdots
n_r)}{n_1^{s_1}\cdots n_r^{s_r}} = \left( \sum_{n_1=1}^{\infty}
\frac{f(n_1)}{n_1^{s_1}} \right) \cdots \left( \sum_{n_r=1}^{\infty}
\frac{f(n_r)}{n_r^{s_r}} \right)
\end{equation*}
\begin{equation} \label{Dir_ser}
\times \prod_p \left(1+ g(p)h(p) \sum_{j=2}^r (-1)^{j-1} f(p^{j-2}) \sum_{1\le i_1< \ldots <
i_j\le r} \frac1{p^{s_{i_1}+\ldots +s_{i_j}}}\right).
\end{equation}

To obtain this result use \eqref{F_homogen} for $k=2$, \eqref{ident_h_2} and \eqref{Dir_ser_g}.
Since the function $(n_1,\ldots,n_r) \mapsto
f(n_1\cdots n_r)$ is multiplicative, we deduce
\begin{equation*}
D:= \sum_{n_1,\ldots,n_r=1}^{\infty} \frac{f(n_1\cdots
n_r)}{n_1^{s_1}\cdots n_r^{s_r}} =\prod_p
\sum_{\nu_1,\ldots,\nu_r=0}^{\infty}
\frac{f(p^{\nu_1+\ldots+\nu_r})}{p^{\nu_1s_1+\ldots+\nu_rs_r}}
\end{equation*}
\begin{equation*}
= \prod_p \sum_{\nu_1,\ldots,\nu_r=0}^{\infty}
\frac{x^{\nu_1+\ldots+\nu_r+1}- y^{\nu_1+\ldots+\nu_r+1}}
{(x-y)p^{\nu_1s_1+\ldots+\nu_rs_r}},
\end{equation*}
where $x=g(p)$ and $y=h(p)$, for short.

Therefore,
\begin{equation*}
D= \prod_p \frac1{x-y} \left(x
\sum_{\nu_1,\ldots,\nu_r=0}^{\infty}
\left(\frac{x}{p^{s_1}}\right)^{\nu_1}\cdots
\left(\frac{x}{p^{s_r}}\right)^{\nu_r} -  y
\sum_{\nu_1,\ldots,\nu_r=0}^{\infty}
\left(\frac{y}{p^{s_1}}\right)^{\nu_1}\cdots
\left(\frac{y}{p^{s_r}}\right)^{\nu_r}\right)
\end{equation*}
\begin{equation*}
= \prod_p \frac1{x-y} \left( x
\left(1-\frac{x}{p^{s_1}}\right)^{-1}\cdots
\left(1-\frac{x}{p^{s_r}}\right)^{-1} - y
\left(1-\frac{y}{p^{s_1}}\right)^{-1}\cdots
\left(1-\frac{y}{p^{s_r}}\right)^{-1}\right)
\end{equation*}
\begin{equation*}
= \prod_p \frac1{x-y} \left(1-\frac{x}{p^{s_1}}\right)^{-1}
\left(1-\frac{y}{p^{s_1}}\right)^{-1} \cdots
\left(1-\frac{x}{p^{s_r}}\right)^{-1}
\left(1-\frac{y}{p^{s_r}}\right)^{-1}
\end{equation*}
\begin{equation*}
\times \left(x \left(1-\frac{y}{p^{s_1}}\right)\cdots
\left(1-\frac{y}{p^{s_r}}\right)- y
\left(1-\frac{x}{p^{s_1}}\right)\cdots
\left(1-\frac{x}{p^{s_r}}\right) \right)
\end{equation*}
\begin{equation*}
= D(g,s_1)D(h,s_1)\cdots D(g,s_r)D(h,s_r)
\end{equation*}
\begin{equation*}
\times \prod_p \left(1-xy \sum_{1\le i<j\le r}
\frac1{p^{s_i+s_j}} +xy \left(x+y\right) \sum_{1\le i<j<k\le r}
\frac1{p^{s_i+s_j+s_k}} - \ldots \right.
\end{equation*}
\begin{equation*}
\left. + (-1)^{r-1}xy\left(x^{r-2}+ x^{r-3}y+\ldots
+xy^{r-3}+y^{r-2}\right) \frac1{p^{s_1+\ldots +s_r}} \right),
\end{equation*}
simplifying by $x-y$, which shows that the computations are valid also
in the case when $x-y=0$, that is, $g(p)=h(p)$ for some primes $p$. This gives
\eqref{Dir_ser}. Formula \eqref{Busche_Raman_gener_inv} is obtained
by expressing the function $(n_1,\ldots,n_r)\mapsto f(n_1)\cdots
f(n_r)$ from the convolutional identity \eqref{Busche_Raman_gener}.
The proof is finished.

Note that for the divisor function $d$ formula \eqref{Dir_ser} gives
\begin{equation*}
\sum_{n_1,\ldots,n_r=1}^{\infty} \frac{d(n_1\cdots
n_r)}{n_1^{s_1}\cdots n_r^{s_r}}
\end{equation*}
\begin{equation} \label{Dir_ser_d}
= \zeta^2(s_1)\cdots \zeta^2(s_r)
\prod_p \left(1+ \sum_{j=2}^r (-1)^{j-1} (j-1) \sum_{1\le i_1<
\ldots < i_j\le r} \frac1{p^{s_{i_1}+\ldots +s_{i_j}}}\right),
\end{equation}
which reduces in the case $r=2$ to
\begin{equation*} \sum_{n_1,n_2=1}^{\infty} \frac{d(n_1n_2)}{n_1^{s_1}n_2^{s_2}} =
\frac{\zeta^2(s_1)\zeta^2(s_2)}{\zeta(s_1+s_2)}.
\end{equation*}

The identity \eqref{Dir_ser_d} corresponding to the cases $r=2$ and
$r=3$ was pointed out in \cite[Sect.\ 6]{KurOch2009}. We also remark
that for the $\sigma$ function the common analytic version of
\eqref{Busche_Raman_i} and \eqref{Busche_Raman_ii} is the formula
\begin{equation*} \sum_{n_1,n_2=1}^{\infty} \frac{\sigma(n_1n_2)}{n_1^{s_1}n_2^{s_2}} =
\frac{\zeta(s_1)\zeta(s_1-1)\zeta(s_2)\zeta(s_2-1)}{\zeta(s_1+s_2-1)}.
\end{equation*}

%%%%%%%%%%%%%%%%%%%%%%%%%%%%%%%%%%%%%%%%%%%%%%%%%%%%%%%%%%%%%%%%%%%%%%%%%
%%%%%%%%%%%%%%%%%%%%%%%%%%%%%%%%%%%%%%%%%%%%%%%%%%%%%%%%%%%%%%%%%%%%%%%%%

\section{Proof of Theorem \ref{Th_2}}

Using \eqref{F_homogen} and \eqref{ident_h} we deduce
\begin{equation*}
\sum_{n_1,n_2=1}^{\infty} \frac{F(n_1n_2)}{n_1^{s_1}n_2^{s_2}}=
\prod_p \sum_{\nu_1,\nu_2=0}^{\infty}
\frac{F(p^{\nu_1+\nu_2})}{p^{\nu_1s_1+\nu_2s_2}}
\end{equation*}
\begin{equation*}
= \prod_p \left( \sum_{\nu_1,\nu_2=0}^{\infty}
\frac1{p^{\nu_1s_1+\nu_2s_2}} \sum_{i=1}^k x_i^{\nu_1+\nu_2+k-1} \prod_{\substack{j=1\\
j\ne i}}^k (x_i-x_j)^{-1} \right)
\end{equation*}
\begin{equation*}
= \prod_p \left(\sum_{i=1}^k x_i^{k-1} \prod_{\substack{j=1\\
j\ne i}}^k (x_i-x_j)^{-1} \sum_{\nu_1,\nu_2=0}^{\infty}
\left(\frac{x_i}{p^{s_1}}\right)^{\nu_1}
\left(\frac{x_i}{p^{s_2}}\right)^{\nu_2}\right)
\end{equation*}
\begin{equation*}
= \prod_p \left(\sum_{i=1}^k x_i^{k-1}
\left(1-\frac{x_i}{p^{s_1}}\right)^{-1}
\left(1-\frac{x_i}{p^{s_2}}\right)^{-1} \prod_{\substack{j=1\\
j\ne i}}^k (x_i-x_j)^{-1} \right)
\end{equation*}
\begin{equation*}
= \prod_p \left( \prod_{\ell=1}^k
\left(1-\frac{x_{\ell}}{p^{s_1}}\right)^{-1}
\left(1-\frac{x_{\ell}}{p^{s_2}}\right)^{-1} \sum_{i=1}^k x_i^{k-1}
 \prod_{\substack{j=1\\ j\ne i}}^k (x_i-x_j)^{-1} \left(1-\frac{x_j}{p^{s_1}}\right)
\left(1-\frac{x_j}{p^{s_2}}\right) \right)
\end{equation*}
\begin{equation*}
= \left( \sum_{n_1=1}^{\infty} \frac{F(n_1)}{n_1^{s_1}} \right)
\left( \sum_{n_2=1}^{\infty} \frac{F(n_r)}{n_2^{s_2}} \right)
\prod_p Q_k(p^{-s_1},p^{-s_2}),
\end{equation*}
where $Q_k(u,v)$ is the polynomial in $u$ and $v$, given by
\begin{equation*}
Q_k(u,v)= \sum_{i=1}^k x_i^{k-1}
 \prod_{\substack{j=1\\ j\ne i}}^k (x_i-x_j)^{-1} \left(1-x_ju \right)
\left(1-x_jv \right)= \sum_{m,n=0}^{k-1} c_{m,n} u^mv^n.
\end{equation*}

Here the coefficients $c_{m,n}$ ($1\le m,n\le k-1$) are given by
\begin{equation} \label{coeff}
c_{m,n} = (-1)^{m+n} \sum_{i=1}^k x_i^{k-1}
 \prod_{\substack{j=1\\ j\ne i}}^k (x_i-x_j)^{-1}
 e^{(i)}_m(x_1,\ldots, x_k) e^{(i)}_n(x_1,\ldots, x_k)
\end{equation}
where $e^{(i)}_m(x_1,\ldots, x_k) =e_m(x_1,\ldots,x_{i-1},x_{i+1},
\ldots,x_k)$ ($1\le i\le k$). We will show that
\begin{equation*}
c_{m,n} =\begin{cases} 1, & m=n=0,\\
(-1)^{m+n-1}e_{m+n}(x_1,\ldots,x_k), & m,n\ge 1, m+n\le k \\ 0, &
\text{otherwise}.
\end{cases}
\end{equation*}

To this end, note that $e^{(i)}_m(x_1,\ldots, x_k)= e_m(x_1,\ldots,
x_k)-x_i e^{(i)}_{m-1}(x_1,\ldots, x_k)$ ($1\le m\le k$), which
leads to the identity
\begin{equation*}
e^{(i)}_m(x_1,\ldots, x_k)= \sum_{\ell=0}^m (-1)^{\ell} x_i^{\ell}
e_{m-\ell} (x_1,\ldots, x_k).
\end{equation*}

Therefore, from \eqref{coeff} we deduce
\begin{equation*}
c_{m,n} = (-1)^{m+n} \sum_{i=1}^k x_i^{k-1}
 \prod_{\substack{j=1\\ j\ne i}}^k (x_i-x_j)^{-1}
\sum_{\ell=0}^m (-1)^{\ell} x_i^{\ell} e_{m-\ell} (x_1,\ldots, x_k)
\sum_{s=0}^n (-1)^s x_i^s e_{n-s} (x_1,\ldots, x_k)
\end{equation*}
\begin{equation*}
= (-1)^{m+n} \sum_{\ell=0}^m \sum_{s=0}^n (-1)^{\ell+s}
e_{m-\ell}(x_1,\ldots, x_k) e_{n-s} (x_1,\ldots, x_k) \sum_{i=1}^k
x_i^{\ell+s+k-1} \prod_{\substack{j=1\\ j\ne i}}^k (x_i-x_j)^{-1}
\end{equation*}
\begin{equation} \label{comput}
= (-1)^{m+n} \sum_{\ell=0}^m (-1)^{\ell} e_{m-\ell}(x_1,\ldots, x_k)
\sum_{s=0}^n (-1)^s e_{n-s}(x_1,\ldots, x_k)
h_{\ell+s}(x_1,\ldots,x_k),
\end{equation}
using again \eqref{ident_h}. For $n=0$ this gives, by
\eqref{rel_e_h},
\begin{equation*}
c_{m,0}= (-1)^{m} \sum_{\ell=0}^m (-1)^{\ell} e_{m-\ell}(x_1,\ldots,
x_k) h_{\ell}(x_1,\ldots,x_k)= \begin{cases} 1, & m=0,\\ 0, & m\ge
1.
\end{cases}
\end{equation*}

We deduce in the same way that $c_{0,n}=0$ for $n\ge 1$. Now let
$m,n\ge 1$. Then the inner sum in \eqref{comput} is, by denoting
$j=\ell+s$,
\begin{equation*}
\sum_{j=\ell}^{n+\ell} (-1)^{j-\ell} e_{n+\ell-j}(x_1,\ldots, x_k)
h_j(x_1,\ldots,x_k)
\end{equation*}
\begin{equation*}
= (-1)^{\ell} \left(\sum_{j=0}^{n+\ell} (-1)^j
e_{n+\ell-j}(x_1,\ldots, x_k) h_j(x_1,\ldots,x_k) \right.
\end{equation*}
\begin{equation*}
\left. - \sum_{j=0}^{\ell-1} (-1)^j e_{n+\ell-j}(x_1,\ldots, x_k)
h_j(x_1,\ldots,x_k)\right)
\end{equation*}
\begin{equation} \label{comput_2}
= (-1)^{\ell-1} \sum_{j=0}^{\ell-1} (-1)^j e_{n+\ell-j}(x_1,\ldots,
x_k) h_j(x_1,\ldots,x_k),
\end{equation}
since the first sum is zero, according to \eqref{rel_e_h}, where
$n+\ell \ge 1+\ell\ge 1$ for every $\ell \ge 0$. For $\ell=0$
\eqref{comput_2} is zero (empty sum). We obtain
\begin{equation*}
c_{m,n}=(-1)^{m+n-1} \sum_{\ell=1}^m e_{m-\ell}(x_1,\ldots,x_k)
\sum_{j=0}^{\ell-1} (-1)^j e_{n+\ell-j}(x_1,\ldots, x_k)
h_j(x_1,\ldots,x_k)
\end{equation*}
and regrouping the terms according to the values $t=\ell-j$,
\begin{equation*}
c_{m,n}=(-1)^{m+n-1} \sum_{t=1}^m e_{n+t}(x_1,\ldots,x_k)
\sum_{j=0}^{m-t} (-1)^j e_{m-t-j}(x_1,\ldots, x_k)
h_j(x_1,\ldots,x_k),
\end{equation*}
where the inner sum is $0$ for $t<m$ and it is $1$ for $t=m$.
Therefore,
\begin{equation*}
c_{m,n}=(-1)^{m+n-1}e_{m+n}(x_1,\ldots, x_k),
\end{equation*}
which is zero for $m+n>k$. This finishes the proof of
\eqref{Busche_Raman_gener_other}. Now
\eqref{Busche_Raman_gener_other_inv} is obtained by expressing the
function $(n_1,n_2)\mapsto F(n_1)F(n_2)$.

%%%%%%%%%%%%%%%%%%%%%%%%%%%%%%%%%%%%%%%%%%%%%%%%%%%%%%%%%%%%%%%%%%%%%%%%%
%%%%%%%%%%%%%%%%%%%%%%%%%%%%%%%%%%%%%%%%%%%%%%%%%%%%%%%%%%%%%%%%%%%%%%%%%

\section*{Acknowledgement} The author gratefully acknowledges support
from the Austrian Science Fund (FWF) under the project Nr.
M1376-N18.

%%%%%%%%%%%%%%%%%%%%%%%%%%%%%%%%%%%%%%%%%%%%%%%%%%%%%%%%%%%%%%%%%%%%%%%%%
%%%%%%%%%%%%%%%%%%%%%%%%%%%%%%%%%%%%%%%%%%%%%%%%%%%%%%%%%%%%%%%%%%%%%%%%%

\medskip

\noindent L\'aszl\'o T\'oth \\
Institute of Mathematics, Universit\"at f\"ur Bodenkultur \\
Gregor Mendel-Stra{\ss}e 33, A-1180 Vienna, Austria \\ and \\
Department of Mathematics, University of P\'ecs \\ Ifj\'us\'ag u. 6,
H-7624 P\'ecs, Hungary \\ E-mail: ltoth@gamma.ttk.pte.hu

\end{document}